\documentclass{amsart}
\usepackage{amssymb}
\usepackage[final]{lpretex}
\usepackage{title}
\usepackage{a4}
\input xy
\xyoption{all}

\EndOfDump

\def\DHpartformat#1{(#1) }

\LBsetstyleof{partno}{arabic}

\DocVersion[Changed to LaTeX]{1}

\let\define\def

  \def\C {{\mathbb C}}

\def\GG {{\mathbb G}}   
  
  \def\P {{\mathbb P}} 
\def\Q {{\mathbb Q}} \def\R {{\mathbb R}}

\def\Z {{\mathbb Z}} 

\define \n {\mathbb N}
\define \z {\mathbb Z}
\define \q {\mathbb Q}
\define \PP {\mathbb P}

\def\sA {{\Cal A}}  \def\sC {{\Cal C}}
 \def\sE {{\Cal E}} \def\sF {{\Cal F}}
  
  \def\sL {{\Cal L}}
 \def\sN {{\Cal N}} \def\sO {{\Cal O}}
 
  \def\sU {{\Cal U}}
\def\sV {{\Cal V}}

\define \cN {\Cal N}
\define \cf {\Cal F}
\define \cg {\Cal G}
\define \cE {\Cal E}
\define \ce {\Cal E}
\define \cc {\Cal C}
\define \cV {\Cal V}
\define \cA {\Cal A}
\define \cK {\Cal K}
\define \cO {\Cal O}
\define \cF {\Cal F}
\define \cn {\Cal N}
\define \cI {\Cal I}
\define \sP {\Cal P}
\def\tA {\widetilde{\Cal A}}

  \def\g {\gamma}  
\def\eps {\epsilon}

\define \x {\xi}
\define \y {\eta}
\define \G {\Gamma}
\define \r {\rho}
\define \w {\omega}

\def \trho {\widetilde {\rho}}

\def \tp {\widetilde{\mathbb P}}
\define \tH {\widetilde H}
\define \tG {\widetilde{\Gamma}}
\define \tW {\widetilde W}
\define \tF {\widetilde F}
\define \tm {\widetilde m}
\define \St {\widetilde S}
\define \Xt {\widetilde X}
\define \tS {\widetilde S}
\define \tpsi {\widetilde \psi}
\define \tL {\widetilde L}
\define \tE {\widetilde E}
\define \tl {\widetilde l}
\define \tA {\widetilde A}
\define \tom {\widetilde\omega}
\define \tT {\widetilde T}
\define \tB {\widetilde B}
\define \tf {\widetilde f}
\define \tsA {\widetilde{\sA}}

\define \tM {\widetilde M}
\define \tphi {\widetilde{\phi}}
\define \trho {\widetilde{\rho}}
\define \tR {\widetilde R}
\define \tp {\widetilde p}
\define \tq {\widetilde q}
\define \tc {\widetilde c}
\define \tsF {\widetilde {\sF}}
\define \tsN {\widetilde {\sN}}
\define \tsU {\widetilde {\sU}}
\define \th {\widetilde h}
\define \tfra{\widetilde{\frak a}}

\def\pd {\partial}

\def \Dx1 {\frac{\pd}{{\pd} x_1}}
\def \Dy1 {\frac{\pd}{{\pd} y_1}}
\def \Dz1 {\frac{\pd}{{\pd} z_1}}
\def \Dx2 {\frac{\pd}{{\pd} x_2}}
\def \Dy2 {\frac{\pd}{{\pd} y_2}}
\def \Dz2 {\frac{\pd}{{\pd} z_2}}

\def\q {\quad} 

\def\from{\longleftarrow}

\def\mapdiagr#1{\Big\searrow\rlap{$\raise 5pt\vbox{{\hbox{$\mkern -15mu\scriptstyle#1$}}}$}}   

\def\mapdiagl#1{\llap{$\raise 5pt\vbox{{\hbox{$\scriptstyle#1\mkern
-15mu$}}}$}\Big\swarrow}              

\def\Mapdiagr#1{\nearrow\rlap{$\lower 5pt\vbox{{\hbox{$\mkern
-15mu\scriptstyle#1$}}}$}} 

\def\Mapdiagl#1{\llap{$\lower 5pt\vbox{{\hbox{$\scriptstyle#1\mkern
-15mu$}}}$}\searrow} 

\def\Mapswr#1{\swarrow\rlap{$\lower 5pt\vbox{{\hbox{$\mkern
-15mu\scriptstyle#1$}}}$}}              

\def\Mapnwl#1{\nwarrow\rlap{$\lower 5pt\vbox{{\hbox{$\mkern
-15mu\scriptstyle#1$}}}$}}

\def \inj {\hookrightarrow}

\define \Rhook {\hookrightarrow}

\def \half {\raise1pt\hbox{$\scriptstyle
        \frac{1}{2}\displaystyle$}}

\def \x{{\sl X}\llap{$\mkern -2mu {\scriptstyle -}$}}
\def \Symm {\operatorname{Sym}}
\def \Res {\operatorname{Res}}

\define \Kod {\operatorname{Kod}}
\define \dimension {\operatorname{dim}}
\define \codim {\operatorname{codim}}
\define \contr {\operatorname{contr}}
\define \rk {\operatorname{rank}}
\define \im {\operatorname{im}}
\define \Mor {\operatorname{Mor}}
\define \Cl {\operatorname{Cl}}
\define \Hilb {\operatorname{Hilb}}
\define \degree {\operatorname{deg}}
\define \mult {\operatorname{mult}}
\define \Aut {\operatorname{Aut}}
\define \NS {\operatorname{NS}}
\define \Gal {\operatorname{Gal}}
\define \ch {\operatorname{char}}
\define \Jac {\operatorname{Jac}}
\define \Km {\operatorname{Km}}
\define \Sec {\operatorname{Sec}}
\define \Stab {\operatorname{Stab}}
\define \Br {\operatorname{Br}}
\define \inv {\operatorname{inv}}
\define \tr {\operatorname{tr}}
\define \Frob {\operatorname{Frob}}
\define \Symn {\operatorname{Sym}^n}
\define \Ev {\sE^\vee}
\define \ordp {\operatorname{ord}_p}
\define \Supp {\operatorname{Supp}}
\define \Ann {\operatorname{Ann}}
\define \disc {\operatorname{disc}}
\define \Lie {\operatorname{Lie}}
\define \embdim {\operatorname{embdim}}

\def\te{\tilde{\eta}}

\def\hod#1#2#3#4{\ensuremath{\if#30 H^{#2}({#1},{\cal O}_{#1}) \else 
 H^{#2}(#1,\Omega^{#3}\if\relax{#4}\relax_{#1}\else _{#1/#4}\fi)\fi}}

\begin{document}

\title{The non-existence of stable Schottky forms}

\author{G.Codogni}
\address{D.P.M.M.S.\\
University of Cambridge\\
Cambridge CB3 0WB U.K.}

\author{N. I. Shepherd-Barron}
\address{Dept. of Mathematics\\
King's College London\\
Strand\\
London WC2R 2LS U.K.}
\email{g.codogni@dpmms.cam.ac.uk, nisb@dpmms.cam.ac.uk}

\maketitle
\let\thefootnote\relax\footnote{To appear in Compositio Mathematica}

\keywordsname{:\ curve, abelian variety, moduli, Schottky problem}


Abstract: We show that there is no stable Siegel modular form
that vanishes on every moduli space of curves.

\begin{section}{Introduction}
Denote by $M_g$ and $A_g$ the coarse moduli spaces of genus $g$ curves
and principally polarized abelian $g$-folds, respectively, and
by $j_g:M_g\to A_g$ the Jacobian 
morphism, that sends every curve
to its Jacobian. This separates geometric points,
according to (a crude version of) the Torelli theorem, so that
we can be somewhat careless in distinguishing between $M_g$
and the Jacobian locus (the image of $M_g$ under $j_g$).

The Satake compactification $A_g^S$ of $A_g$
is the minimal complete normal variety that contains $A_g$ as an 
open subvariety
and is also characterized by the property that the $\Q$-line bundle
$\sL$ (sometimes denoted by $\omega$; however, this notation conflicts
with the use of $\omega$ to denote the canonical bundle)
of weight $1$ (Siegel) modular forms on $A_g$ extends
to an ample $\Q$-line bundle on $A_g^S$. There is a stratification
$A_g^S=\cup_{0\le h\le g}A_h$, so that the boundary 
$\partial A_g^S=A_g^S-A_g$ of $A_g^S$
is $\partial A_g^S=\cup_{0\le h\le g-1}A_h$. It follows that
$A_{g-1}^S$ is then the normalization of $\partial A_g^S$;
in fact, $A_{g-1}^S=\partial A_g^S$, from the surjectivity
of the Siegel operator $\Phi$ (some of whose basic properties
recalled below) for modular forms of
sufficiently high weight. In this way we can regard $A_g^S$
as a subvariety of
$A_{g+m}^S$ for every positive $g,m$.

Define the Satake compactification $M_g^S$ of $M_g$
to be the closure of the Jacobian locus in $A_g^S$. So 
[H], [M] the
geometric points of $A_g^S$ correspond to ppav's of dimension
at most $g$ and those of $M_g^S$ correspond to products of Jacobians,
where the genera of the curves in question sum to at most $g$.
So the intersection $M_{g+m}^S\cap A_g^S$, taken inside
$A_{g+m}^S$ is exactly $M_g^S$, as sets. The main result
of this paper (which was inspired by a recent result of
Grushevsky and Salvati Manni [G-SM] that we recall below) 
is that this intersection is far from being
transverse, however.

\begin{theorem}\label{non-transverse} 
$M_{g+m}^S\cap A_g^S$ contains
the $m$th order infinitesimal neighbourhood of
$M_g^S$ inside $A_g^S$. 
\noproof
\end{theorem}

By definition, if an integral subscheme $X$ of a scheme $Y$
is defined by an ideal $I$, then the $m$th order
infinitesimal neighbourhood of
$X$ in $Y$ is the subscheme of $Y$ defined by the 
$(m+1)$st symbolic power $I^{[m+1]}$ of $I$.
When $X$ and $Y$ are both regular,
then $I^{[m+1]}=I^{m+1}$,
the ordinary $(m+1)$st power.
If, moreover, $X$ and $Y$ are smooth over a field of characteristic
zero, if $x_1,...,x_n$ are holomorphic (or formal, or
{\'e}tale) co-ordinates on $Y$ such that $I$ is generated
by $x_1,...,x_r$, then $I^m$ consists of those
functions $f$ on $Y$ such that $f$ and all its partial derivatives
with respect to $x_1,...,x_r$ up to and including
those of order $m-1$ vanish along $X$.

\begin{corollary} If $g$ is fixed, then
$\cup_{m\ge 0}(M^S_{g+m}\cap A_g^S)$ equals the
formal completion $(A_g^S)\widehat{}$
of $A_g^S$ along $M_g^S$.
\noproof
\end{corollary}

Recall that the Schottky problem, 
in its general form, is the problem of distinguishing
Jacobians from other principally polarized abelian varieties. One
classical approach is to seek \emph{Schottky forms},
that is, scalar-valued Siegel modular forms on $A_g$ that
vanish on the Jacobian locus, or, equivalently, 
forms on $A_g^S$ that vanish on $M_g^S$.

The normalization $\nu:A_{g}^S\to\partial A_{g+1}^S$
gives a restriction map $\Phi$, which coincides with 
the \emph{Siegel operator}, from the vector space 
$[\G_{g+1},k]=H^0(A_{g+1}^S,\sL^{\otimes k})$ of
weight $k$ forms on $A_{g+1}$ (equivalently, on $A_{g+1}^S$)
to $[\G_{g},k]$. This is surjective if $k$ is even and $k>2g$
(\cite{Fr2}, p. $64$)
so that, since $\sL$ is ample on the Satake compactification,
$\nu$ is an isomorphism. 

In terms of holomorphic functions
on the Siegel upper half-planes $\frak H_{g+1}$ and $\frak H_g$
of degrees $g+1$ and $g$ respectively,
$\Phi$ is defined by
$$\Phi(F)(\tau)=\lim_{t\to+\infty}F(\tau\oplus it),$$
where the direct sum of two square matrices has its obvious meaning.
In terms of a Fourier expansion
$$F(T)=\sum_S a(S)\exp\pi i\tr(ST),$$
where $T\in\frak H_{g+1}$ and $S$ runs over the positive semi-definite
symmetric integral matrices with even diagonal,
$\Phi$ is given by
$$\Phi(F)(\tau)=\sum_{S_1} a(S_1\oplus 0)\exp \pi i\tr(S_1\tau).$$
Or $\Phi$ can be calculated by first restricting to
a copy of $\frak H_g\times\frak H_1$ in $\frak H_{g+1}$
that arises as a cover of the image of $A_g\times A_1$ in
$A_{g+1}$ and then letting $\tau\in\frak H_1$ tend to $i\infty$.

In genus $4$ Schottky found one such Schottky form explicitly [S]; 
later Igusa [I1] showed that, if 
$$F_g=\theta_{E_8^2,g}-\theta_{D_{16}^+,g},$$
the difference of the theta series 
in genus $g$ associated to the two distinct positive even
unimodular quadratic forms $E_8^2$ and $D_{16}^+$ of rank
$16$, then the form discovered by Schottky is an
explicit rational multiple of $F_4$. He also [I2]
showed that
it is reduced and irreducible, and so cuts out
exactly the Jacobian locus in $A_4$.

(Recall that if a positive even unimodular quadratic form $Q$
of rank $k$ is regarded as a lattice $\Lambda_Q$ in Euclidean
space $\R^k$, then
the theta series $\theta_{Q,g}$ is defined by
$$\theta_{Q,g}(\tau)=\sum_{x_1,...,x_g\in\Lambda_Q}
\exp\pi i\sum_{p,q=1}^gQ(x_p,x_q)\tau_{pq}.$$
This lies in the space $[\Gamma_g,k/2]$ 
and satisfies the formula
$\Phi(\theta_{Q,g+1})=\theta_{Q,g}$.)

Grushevsky and Salvati Manni \cite{G-SM} have shown that 
the genus $5$ Schottky form $F_5$ does not vanish along $M_5$.
(They go on to prove that $F_5$ cuts out exactly
the trigonal locus in $M_5$.)
They did this by proving that if $F_5$ did vanish along
$M_5$, then $F_4$ would vanish with multiplicity at least
$2$ along $M_4$ (that is, $F_4$ and all its first partial
derivatives would vanish), which would contradict
Igusa's result on reducedness. Theorem \ref{non-transverse}
is a generalization of this.

To begin with,
say that a scalar-valued Siegel modular
form $F$ on $A_g$ vanishes with multiplicity at least $m$ along
$M_g$ if $F$ and all its partial derivatives
with respect to the co-ordinates $\tau_{pq}$ on $\frak H_g$
of order at most $m-1$ vanish along $M_g$ (rather,
along the inverse image of $M_g$ in $\frak H_g$).

(At this point there is a slight conflict between 
the language of stacks and that of varieties.
On the one hand
the Satake compactification,
which is fundamental here, is a coarse object and on the other
the variables $\tau_{pq}$ are local co-ordinates
on the stack $\sA_g$ but not on the coarse space $A_g$.)

\begin{theorem}\label{vanishing} Suppose that $F=F_{g+1}$
is a scalar-valued Siegel modular form on $A_{g+1}$ that
vanishes with multiplicity at least $m\ge 1$
along the Jacobian locus $M_{g+1}$ in $A_{g+1}$. Then
$F_g=\Phi(F_{g+1})$ vanishes with multiplicity at 
least $m+1$ along $M_g$.
\noproof
\end{theorem} 

Note that Theorem \ref{non-transverse} follows at once from
Theorem \ref{vanishing} and a lemma in commutative algebra,
Lemma \ref{primary} below.

Finally, we restate this in terms of Freitag's description \cite{Fr1}
of the ring of stable Siegel modular forms. He showed that, for a
fixed even integer $k$, the Siegel map
$\Phi:[\G_{g+1},k]\to [\G_g,k]$ is an isomorphism for all
$g>2k$. That is, the vector spaces $[\G_g,k]$ stabilize
to a vector space $[\G_\infty,k]$ as $g$ increases,
and is the space of sections of a $\Q$-line bundle
$\sL^{\otimes k}$ on $A_\infty^S$. Let
$\sN$ denote the restriction of $\sL$ to $M_g^S$.

Put $A(\G_g)=\oplus_k[\G_g,k]$,
the graded ring of Siegel modular forms on the moduli
space $A_g$ or on the Satake compactification
$A_g^S$. Then $A=\oplus_k [\G_\infty,k]$ is
is an inverse limit, in the category of graded rings:
$$A=\lim_{\stackrel{\from}{g}}A(\G_g).$$
Freitag proved also that $A$ is the polynomial ring over
$\C$ on the set of theta series $\theta_Q$, where $Q$
runs over the set of equivalence classes of
indecomposable even, positive and
unimodular quadratic forms over $\Z$.

Define the stable Satake compactification
$A_\infty^S$ by 
$$A_\infty^S=\cup_g A_g^S
=\lim_{\stackrel{\to}{g}}A_g^S.$$
Define $M_\infty^S=\cup_g M_g^S$,
similarly.

\begin{corollary}
The homomorphism from $A$ to the graded ring
$\oplus_k H^0(M_\infty^S,\sN^{\otimes k})$ that is induced
by the inclusion $M_\infty^S\inj A_\infty^S$
is injective.
That is, there are no stable Schottky forms.
\begin{proof}
An element $F$ of the kernel would restrict, in each genus
$g$, to a scalar-valued modular form $F_g$ on $A_g$ that vanished to
arbitrarily high order along $M_g$. Then $F_g=0$ for all $g$,
and then $F=0$.
\end{proof}
\end{corollary}
\begin{corollary}\label{genus} If $P,Q$ are positive even
unimodular quadratic forms that are not equivalent,
then there exists a curve $C$ whose period matrix distinguishes
between them.\noproof
\end{corollary}
If $P,Q$ have rank $g$, then there is a period
matrix $\tau$ in $\frak H_g$
such that $\theta_{P,k}(\tau)\ne\theta_{Q,k}(\tau)$. However,
it is not clear how to find the genus of the curve whose
existence is given by Corollary \ref{genus},
let alone how to identify a particular such curve.
However, more recently we have [C], [SB] shown
that $C$ can be taken to be a general point of
the trigonal locus (of some genus, as yet undetermined), 
but not the hyperelliptic locus.
That is, there exist stable modular forms
vanishing on every hyperelliptic locus,
but there are no such forms that vanish
on every trigonal locus.
\end{section}
\bigskip

\begin{section}{Fay's degenerating families}
\medskip
We make a slight extension of a
construction by Fay \cite{F}, pp. 50-54, of $1$-parameter families
of genus $g+1$ curves that degenerate to an irreducible nodal
curve of geometric genus $g$. His construction 
includes a calculation of the
period matrix of the general member of the family,
modulo the square of the parameter. 
The formula that he gives is, however, mistaken,
as was pointed out by Yamada [Y], who gave the correct
version. The 
reason for the extension is to permit a 
rescaling of local co-ordinates by
non-zero parameters $\lambda,\mu$;
Fay's construction, with its original wording, only
permits this rescaling when $|\lambda|=|\mu|=1$.

Start with a curve $C$ of genus
$g$. Let $\sV$ be the infinite-dimensional variety
whose points are quadruples
$(a,b,z_a,z_b)$, where $a,b$ are distinct
points on $C$ and $z_a,z_b$ are local
holomorphic co-ordinates on $C$ at $a,b$ respectively.

The $2$-torus $\GG_m^2$ acts on $\sV$ by 
$$(\lambda,\mu)(a,b,z_a,z_b)=(a,b,\lambda^{-1}z_a,\mu^{-1}z_b).$$
Fix a non-empty finite-dimensional (in order to avoid
irrelevant difficulties) and smooth subvariety $V$
of $\sV$ that is preserved under this 
torus action and that maps onto the complement $U$
of the diagonal in $C\times C$.

We want to construct a family of morphisms
$\{f_v:\sC_v\to \Delta\}_{v\in V}$ that is parametrized
by $V$, where $\Delta$
is a complex disc centred at $0$, each $\sC_v$ is a smooth
complex surface, each $f_v$ is proper and each fibre 
over $0$ is
a nodal curve $C/(a\sim b)$, every
other fibre is a smooth curve of genus $g+1$ and the parametrization
is holomorphic in $V$.

It is clearer to run through the construction without
referring to the parameter space $V$. So fix the data
$C,a,b,z_a,z_b$ and choose $\delta>0$ 
such that there are disjoint neighbourhoods $U^a$ of $a$ 
and $U^b$ of $b$ such that $z_a:U^a\to\C$ and $z_b:U^b\to\C$
are each an isomorphism to some open set that 
contains a disc of radius $\delta$ centred at $z_a(a)=0$
and $z_b(b)=0$, respectively.

Let $\Delta_\delta$, $D_{\delta^2}$ denote
complex discs of radius $\delta$, ${\delta^2}$, respectively.

Take $W=W_\delta$ to be the open subset of $C\times D_{\delta^2}$ 
obtained by deleting the two closed subsets
$$\{(p,t)\big\vert p\in U^a,\ 0\le\delta |z_a(p)|\le|t|\le\delta^2\},$$
$$\{(q,t)\big\vert q\in U^b,\ 0\le\delta |z_b(q)|\le|t|\le\delta^2\}.$$

\begin{lemma}\label{shrink}
If $\eps<\delta$ then $W_\eps\subset W_\delta$.
\noproof
\end{lemma}

In $W_\delta$, define 
open subsets
$$W^a=W^a_\delta=\{(p,t)\big\vert p\in U^a,\ 0<|z_a(p)|<\delta\ \mathrm{and}\ 
|t|< \delta|z_a(p)|\},$$
$$W^b=W^b_\delta=\{(q,t)\big\vert q\in U^b,\ 0<|z_b(q)|<\delta\ \mathrm{and}\ 
|t|< \delta|z_b(q)|\}.$$
Consider the complex surface 
$S=S_\delta\subset (\Delta_\delta)^2\times D_{\delta^2}$
defined by the equation $XY=t$, where $X,Y$
are co-ordinates on the two copies of $\Delta_\delta$
and $t$ is a co-ordinate on $D_{\delta^2}$.
Then there are isomorphisms
$$W^a_\delta\to S-(X=0): (p,t)\mapsto (z_a(p),t/z_a(p),t),$$
$$W^b_\delta\to S-(Y=0): (q,t)\mapsto (t/z_b(q),z_b(q),t).$$
Together these define an {\'e}tale morphism
$j:W^a_\delta\cup W^b_\delta\to S$, where
the union is the disjoint union, taken inside
$C\times D_{\delta^2}$. Let
$i:W^a_\delta\cup W^b_\delta\to W_\delta$ be the inclusion.

If $Z$ is a subspace of a space $X$, then $\overline{Z}$
denotes the closure of $Z$ in $X$.

\begin{lemma} $(i,j):W^a_\delta\cup W^b_\delta\to W\times S$
is a closed embedding.
\begin{proof} It is enough to show that the image of
$W^a_\delta$ in ${\overline{W^a_\delta}}\times S$ is closed.
Now points in ${\overline{W^a_\delta}}\times S$ are of the form
$(p,t_1,X,Y,t_2)$ with
\begin{eqnarray*}
\delta\ge |z_a(p)|\ge & t_1/\delta,\ t_2=t_1,\ X=z_a(p),\ Y=t_2/z_a(p),\\
& |X|,|Y|<\delta,\ |t_2|\le \delta^2,\ XY=t_2.
\end{eqnarray*}
But these conditions force 
$\delta>|z_a(p)|=|t_2|/|Y|> |t_2|/\delta,$
and we are done.
\end{proof}
\end{lemma}

Now define $\sC=\sC_\delta$ by glueing $W_\delta$ to $S_\delta$ by the
inclusion $i$ and the {\'e}tale map $j$.
By the lemma, $\sC$ is Hausdorff
(Bourbaki, Top. G{\'e}n. I.9, p.57, Prop. 4.),
and by construction there is a morphism $f:\sC\to D_{\delta^2}$
whose fibre over $0$ is the nodal curve $C/(a\sim b)$.

\begin{lemma} $f$ is proper.
\begin{proof} It is enough to show that, for any
$r\in(0,\delta)$, the inverse image $Z_r=f^{-1}({\overline{D_{r^2}}})$
is compact.
By construction, $Z_r$ is the union of the two
compact spaces ${\overline{W^1_\delta}}$ and
${\overline{S_r}}$, where the subset
${\overline{S_r}}$ of $S_\delta$ is defined by
$|X|,|Y|\le r$.
\end{proof}
\end{lemma} 

\begin{lemma} The restriction of $f:\sC_\delta\to D_{\delta^2}$
to the germ of the pair $(D_{\delta^2},0)$ is independent
of $\delta$.
\begin{proof} This follows from the facts that,
by Lemma \ref{shrink} above, $\sC_\eps$ is open in
$\sC_\delta$, and that $C/(a\sim b)$ is proper.
\end{proof}
\end{lemma}

Note that, by construction,
$W$ is open in $C\times D_{\delta^2}$,
the image of the projection $pr_1:W\to C$
is exactly $C-\{a,b\}$ and
there is an {\'e}tale morphism $\pi:W \to\sC$. 

Given cycles $A_i,B_j$ on $C$ that represent
a symplectic basis of $H_1(C,\Z)$ and are disjoint from
$\{a,b\}$, we can then regard the $A_i,B_j$ as cycles on
$\sC_t$ that represent part of a symplectic
basis of $H_1(\sC_t,\Z)$ for $t\ne 0$ by taking
$pr_1^{-1}(A_i)\cap pr_2^{-1}(t)=A_i\times\{t\}$
and the same thing for $B_j$. Define the cycle $A_{g+1}$
on $\sC_t$ by $A_{g+1}=\partial U^b\times\{t\}$;
then $(A_1,...,A_{g+1},B_1,...,B_g$ can be extended
to a symplectic basis of $H_1(\sC_t,\Z)$ where
$B_{g+1}$ projects to a cycle on the nodal curve
$\sC_0=C/(a\sim b)$ that passes through
the node.

We want to extend this construction of a single
degenerating pencil  
$f:\sC\to D$ of curves to the construction
of a family of such pencils, where the parameter
space is $V$ and the pencil depends
holomorphically on $V$.
This is merely a matter of enhancing the notation
that we have just used, and the details are omitted.
The end result of the construction is a parameter
space $D$ that is an open neighbourhood of
$V\times\{0\}$ in $V\times\C$ and a proper
flat morphism
$\sC\to D$ from an $(n+1)$-dimensional
complex manifold to a complex $n$-manifold that 
is smooth outside $V\times\{0\}$ and whose
restriction to $V\times\{0\}$ is trivial,
with fibre $C/(a\sim b)$.

Now we can follow Fay and Yamada.

We have already chosen $1$-cycles $(A_i,B_j)_{i,j=1,...,g}$
that represent a symplectic basis of $H_1(C,\Z)$.
Take the corresponding normalized basis $(\omega_q)$ 
of $H^0(C,\Omega^1)$. (``Normalized'' means that
$\int_{A_p}\omega_q=\delta_{pq}$ rather than
$2\pi i\delta_{pq}$, which latter is the sense in which
Fay uses the word.) 
Denote by $\tau$ the resulting period
matrix of $C$: that is, $\tau_{pq}=\int_{B_p}\omega_q$.

Also let $\omega_{g+1}=\omega_{b-a}$,
the unique rational $1$-form 
on $C$ whose polar divisor is $a+b$, such that
$\int_{A_p}\omega_{b-a}=0$ for all $p$
and $\Res_b\omega_{b-a}=-\Res_a\omega_{b-a}=\frac{1}{2\pi i}$.

Define scalars $v_p(a)$, etc., by
$v_p(a)=\frac{\omega_p}{dz_a}(a)$; then
for each $p$ the map
$(a,b,z_a,z_b)\mapsto (v_p(a),v_p(b))$
is a holomorphic function $V\to\C^2$.
Also,
take a co-ordinate $t$ on $\C$, so that
$V\times\{0\}$ is the divisor in $D$
defined by $t=0$.

When a curve varies in a holomorphic family, its period
matrix is a holomorphic function of the parameters,
and for the degenerating family just constructed Fay makes 
this explicit, as follows.

\begin{theorem} (Corollary 6 of [Y]) 
After passing to a suitable infinite cyclic
cover of $D-(V\times\{0\})$ there is a symplectic basis of the 
homology of a smooth fibre $\sC_{v,t}$ with respect to
which the period matrix $T=T(v,t)$ of $\sC_{v,t}$
can be written in $2\times 2$ block form 
$$
T=
\left[ {\begin{array}{cc}
{\tau +t\sigma} & {AJ(b)-AJ(a) +ts}\\
{{}^t({AJ(b)-AJ(a) +ts})} & {\frac{1}{2\pi i}(\log t+c_1+c_2t)}\\
\end{array}} \right] +O(t^2).
$$
Here, $AJ$ is the Abel-Jacobi map from the curve
$C$ to its Jacobian,
so that $AJ(b)-AJ(a)$ is the vector
$(\int_a^b\omega_p)$; $s=(s_p)$ is some vector-valued
holomorphic function on $V$ whose explicit form we do not need; 
$c_1,c_2$ are holomorphic functions on $V$ but independent of $t$;
$O(t^2)$ is a holomorphic function on
$\Delta$ that vanishes modulo $t^2$;
and the $g\times g$ matrix $\sigma=(\sigma_{pq})$ is given by
$$\sigma_{pq}=2\pi i\left(v_p(a)v_q(b)+v_q(a)v_p(b)\right).$$
\begin{proof} This is only a matter of verifying
that Fay's calculation goes through in our, 
slightly more general, context. Note, however, that
Fay uses the symbol $v_p$ to denote the normalized
holomorphic $1$-form $\omega_p$, while his expression
$v_p(a)$ must be interpreted as $\frac{\omega_p}{dz_a}(a)$.

The calculation goes as follows. Fix a point in the parameter
space $V$; then we have a degenerating family
$\sC\to D$ of curves, defined locally by an equation
$XY=t$, where $X,Y$ are co-ordinates on the smooth complex surface $\sC$.
By construction,
$$X=z_a(p_a)=t/z_b(p_b),\ Y= t/z_a(p_a) = z_b(p_b).$$
Let $h:C\to\sC$ be the normalization of the degenerate fibre
and $h_a:U^a\to\sC$, $h_b:U^b\to\sC$ be its restriction
to the two given charts on $C$.
Then $h^a$ is defined by $Y=0,2x=X=z_a$ and
$h^b$ by $X=0, 2x=z_b=Y.$

There are holomorphic $2$-forms $\Omega_i$
on the complex surface $\sC$, for $i=1,...,g+1,$
such that if we define
$$u_i(\lambda)=\Res_{\sC_\lambda}\frac{\Omega_i}{t-\lambda},$$
then $(u_i(\lambda))_{i=1,...,g+1}$ is a basis of the space of
holomorphic $1$-forms on $\sC_{\lambda}$, normalized with
respect to the cycles $A_1,...,A_{g+1}$
on $\sC_\lambda$. For $i\le g$, $h^*u_i(0)=\omega_i$,
the normalized $1$-form on $C$. 

Define $W_\lambda$ to be the Riemann surface defined in $W$
by $t-\lambda=0$; this equals $\pi^{-1}(\sC_\lambda)$.
So $W_\lambda$ possesses an {\'e}tale map
$\pi:W_\lambda\to C_\lambda$ and an {\'e}tale map $pr_1:W_\lambda\to C$.
So on $W_\lambda$ there exist $1$-forms 
$\te_i(\lambda)=\pi^*u_i(\lambda)-pr_1^*\omega_i$.
Switch notation from $\lambda$ to $t$.
Note that $\te_i(0)=0$. Now switch notation from $\lambda$ to $t$
and define $\eta_i$ by
$$\eta_i=\lim_{t\to 0}\frac{\te_i}{t}.$$
Then $\eta_i$ is an intrinsic definition of
$\frac{\partial u_i(t)}{\partial t}$.

We can expand $\Omega_i$ in terms of $X,Y$ as
$$\Omega_i=-\phi(X,Y)dX\wedge dY=-\sum c_{m,n}X^mY^ndX\wedge dY.$$
Then
$$u_i(t)= \Res_{\sC_t}\frac{\Omega_i}{XY-t}
=\sum c_{{n+p},n}X^{p-1}t^ndX=\sum c_{{n+p},n}z_a^{p-1}t^ndz_a,$$
where the sum are over $n,p$ with $n\ge 0$, $p\in\Z$ and $n+p\ge 0$.
It follows that $\omega_i=u_i(0)=\sum c_{p,0}z_a^{p-1}dz_a,$
so that $c_{0,0}=0$ and
$\omega_i=\sum_{p\ge 0} c_{p+1,0}z_a^pdz_a$ and then
$$v_i(a)=c_{1,0}.$$

By definition, $\eta_i$ then given by
$$\eta_i=\left(\sum nc_{n+p,n}z_a^{p-1}t^{n-1}dz_a\right)\vert_{t=0},$$
which gives
$$\eta_i=\left(c_{0,1}z_a^{-2}+c_{1,1}z_a^{-1}+\cdots\right)dz_a.$$
Hence $\frac{z_a^2\eta_i}{dz_a}(a)=c_{0,1}$ and $v_i(a)=c_{1,0}$.

An exactly similar calculation, after exchanging $a$ with $b$
and $X$ with $Y$, gives
$\frac{z_b^2\eta_i}{dz_b}(b)=-c_{1,0}$ and
$v_i(b)-c_{0,1}$. (The signs appear because computing residues
in terms of $Y$ instead of $X$ changes the sign.)
That is,
$$\frac{z_a^2\eta_i}{dz_a}(a)=-v_i(b),\
\frac{z_b^2\eta_i}{dz_b}(b)=-v_i(a).$$

The formulae above show that $\eta_i$ is meromorphic.
Differentiating the identity $\int_{A_p}u_i=\delta_{ip}$
gives $\int_{A_p}\eta_i=0$ for all $p\le g$,
and the residues of $\eta_i$ vanish because
the cycle
$A_{g+1}$ equals $\partial U^b\times\{t\}$, by
construction, and the differentials
$u_q$ for $q\le g$ are normalized,
so that $\int_{A_{g+1}}u_q=0$; differentiating
this with respect to $t$, evaluating at $t=0$ and
then pulling back to $C$ gives
$$2\pi i\Res_b\eta_q=\int_{\partial U^b}\eta_q=0.$$
Then $\Res_a\eta_q=0$ also, since the residues
of a meromorphic $1$-form sum to zero.

Now recall the
bilinear relations between holomorphic forms and
those of the second kind
([Spr], p. 260, Theorem 10-8):
if $\phi$ is a meromorphic $1$-form with principal part 
$\frac{\lambda_P}{z_P^2}dz_P$
at each of its poles $P$ (so that, in particular,
$\phi$ has only double poles and all its residues
vanish), where $z_P$ is a local co-ordinate at $P$, 
and if $\omega$ is a holomorphic $1$-form
with $\frac{\omega}{dz_P}(P)=c_P$, then
$$\sum_{j=1}^g\left(\int_{A_j}\omega\int_{B_j}\phi-
\int_{B_j}\omega\int_{A_j}\phi\right)
=2\pi i\sum_P\lambda_P c_P.$$
(Note that the LHS is exactly the cup product
$[\omega]\cup [\phi]$ of the cohomology
classes in $H^1(C,\C)$ defined by these forms,
so the bilinear relations give a formula for the cup
product as a sum of local contributions.)

Take $\omega=\omega_k$ and $\phi=\eta_i$; we get
$$\int_{B_k}\eta_i=-2\pi i\left(v_i(a)v_k(b)+v_k(a)v_i(b)\right).$$
But $\int_{B_k}\eta_i$ is exactly 
the entry $\sigma_{ik}$ of the matrix $\sigma$
appearing in the formula for $T(t)$.

Finally, the entry $T_{g+1,g+1}(t)=\frac{1}{2\pi i}(\log t + c_1 + c_2 t)$
for the reasons of monodromy that Fay gives.
\end{proof}
\end{theorem}

\end{section}
\bigskip
\begin{section}{The failure of transversality}
\medskip
Here we prove Theorem \ref{vanishing}. Recall its statement:

\begin{theorem}\label{3.1} (= Theorem \ref{vanishing})
If $F_{g+1}$ has multiplicity at least $m$ along $M_{g+1}$
then $F_g$ has multiplicity at least $m+1$ along $M_g$.
\begin{proof} Suppose that $N_{g+1}(\{x_{ij}\})$ is a homogeneous
polynomial of degree $d$ in the entries $x_{ij}$ of a symmetric 
$(g+1)\times(g+1)$ matrix $X$. Our hypothesis is that for all $d\le m-1$
and for all such $N_{g+1}$, the partial derivative
$$N_{g+1}(F_{g+1}):=
N_{g+1}\left(\left\{\frac{\partial}{\partial T_{pq}}\right\}\right)(F_{g+1})$$
vanishes along $M_{g+1}$ (rather, its inverse image in 
$\frak H_{g+1}$)
for $T=(T_{pq})\in\frak H_{g+1}$.

Given such $N_{g+1}$, we let $N_g$ denote the polynomial obtained
from it by setting the bottom row and last column of $X$ equal to zero.
Our goal is to show that for every such $N_g$ of degree $m$,
the partial derivative $N_g(F_g)$ vanishes at every point
$\tau$ in $\frak H_g$ that comes from a curve of genus $g$.

For any positive integer $n$, let $S_n$ denote
the set of $n\times n$ integer matrices
that are symmetric, positive semi-definite and 
whose diagonal entries are even. Then
recall that every Siegel modular form $F=F_{g+1}(T)$ 
of degree $g+1$ over a ring $R$
has a Fourier expansion
$$F(T)=\sum_{X\in S_{g+1}}a(X)\exp\pi i\tr(XT)
=\sum_{X\in S_{g+1}}a(X)\exp\pi i\sum_{p,q=1}^{g+1}x_{pq}T_{pq}.$$
We write $X=(x_{pq})$ for
$X\in S_{g+1}$. The Fourier coefficients $a(X)=a_F(X)$ lie in $R$.
For us, $R=\C$.

Take $T$ as above and take $N$ to have degree $m-1$; then 
$$\frac{1}{(\pi i)^{m-1}}N_{g+1}(F_{g+1})(T) =
\sum_{X\in S_{g+1}}a(X)N_{g+1}(\{x_{pq}\})\exp\pi i\sum_{p,q=1}^{g+1}
x_{pq}T_{pq}.$$
Our aim is to examine the coefficient of $t$ in the expansion
of this expression in powers of $t$, so calculate modulo $t^2$.
Since $\exp 2\pi i T_{g+1,g+1}=\g_1\g_2^t t$ modulo $t^2$, where
$\g_j=\exp c_j$, it follows that, modulo $t^2$, we can write
$$\frac{1}{(\pi i)^{m-1}}N_{g+1}(F_{g+1})(T) =\sum_{x_{g+1,g+1}=0}{}+
\sum_{x_{g+1,g+1}=2}{},$$
where $\sum_{x_{g+1,g+1}=2r}{}$ denotes the sum
over matrices $X=(x_{ij})\in S_{g+1}$ whose $(g+1,g+1)$ entry 
$x_{g+1,g+1}$ equals $2r$ (merely because
all terms with $x_{g+1,g+1}\ge 4$ vanish modulo $t^2$).

\begin{lemma} If $X\in S_{g+1}$ and $x_{g+1,g+1}=0$, then the right hand
column and bottom row of $X$ are both zero.
\begin{proof} Immediate consequence of semi-positivity.
\end{proof}
\end{lemma}

Therefore
$$\sum_{x_{g+1,g+1}=0}{}=\sum_{X\in S_g}a(X)N_g(\{x_{pq}\})
\exp\pi i\sum_{p,q=1}^gx_{pq}(\tau_{pq}+t\sigma_{pq})$$
and
\begin{eqnarray*}
\sum_{x_{g+1,g+1}=2}{}= 
& t\g_1\g_2^t \sum_{X\in S_{g+1},x_{g+1,g+1}=2}a(X)N_{g+1}(\{x_{pq}\})\\
&.\left(\exp 2\pi i\sum_{p=1}^gx_{p,g+1}\int_a^b\omega_p\right)
\left(\exp \pi i\sum_{p,q=1}^g x_{pq}\tau_{pq}\right)
\end{eqnarray*}
since we are calculating modulo $t^2$.
So the coefficient of $t$ that we seek is $A+\g_1B,$ where
$$A=\sum_{x_{g+1,g+1}=0}a(X)N_g(\{x_{pq}\})
\left(\pi i\sum_{p,q=1}^gx_{pq}\sigma_{pq}\right)
\left(\exp\pi i\sum_{p,q=1}^gx_{pq}\tau_{pq}\right)$$
and
$$B=\sum_{x_{g+1,g+1}=2}a(X)N_{g+1}(\{x_{pq}\})
\left(\exp 2\pi i\sum_{p=1}^gx_{p,g+1}
\int_a^b\omega_p\right)\left(\exp \pi i\sum_{p,q=1}^g x_{pq}\tau_{pq}\right).$$
The quantities $A,B,\g_1$ are holomorphic
functions on $V$ and, by assumption, $A+\g_1B$ vanishes identically.

Now rescale the local co-ordinates $z_a,z_b$. That is,
start with local co-ordinates $\zeta_a,\zeta_b$ and then
take $z_a=\lambda^{-1}\zeta_a$ and $z_b=\mu^{-1}\zeta_b$.
Such a rescaling will produce a different family $\sC\to \Delta$,
but the quantity $A+\g_1B$ will still vanish
for the rescaled family. Moreover,
$B$ is invariant under this
rescaling, as is revealed by a cursory inspection.
Also $c_1$ is a holomorphic function of $\lambda,\mu$
because the entries of a period matrix are holomorphic
functions of the parameters.

On the other hand, substituting 
$$\sigma_{pq}=-2\pi i\lambda\mu\left(\frac{\omega_p}{d\zeta_a}(a)
\frac{\omega_q}{d\zeta_b}(b)+
\frac{\omega_q}{d\zeta_a}(a)\frac{\omega_p}{d\zeta_b}(b)\right)$$
into the expression above for $A$
shows that $A$ can be written as
$$A=D\lambda\mu,$$
where $D$ is independent of $\lambda,\mu$.
So we have an identity
$$D\lambda\mu=-B\exp(c_1(\lambda,\mu))$$
of holomorphic functions on the $2$-dimensional
torus $\GG_m^2=\Sp\C[\lambda^\pm,\mu^\pm]$,
where we regard $B,D$ as constants
(constant as functions on $\GG_m^2$, that is).

\begin{lemma} Suppose that $f$ is a rational function on 
a complex algebraic variety $X$ and that there is a holomorphic 
function $h$ on some Zariski open subset $U$ of $X$
such that $f=\exp h$ on $U$. Then $f$ is constant.
\begin{proof} It is enough to show that $f$ is constant on a general 
curve in $X$. So we can assume that $\dim X=1$, and then that 
$X$ is a compact Riemann surface. If $f$ is not constant, then it
has a zero, say at $P$, and in some neighbourhood $U$ of $P$ with
a co-ordinate $z$ we have $f=z^nf_1$ with $f_1$ holomorphic
and invertible on $U$, and $n>0$. Then $f_1=\exp h_1$ with $h_1$ holomorphic
on $U$, and $h$ is holomorphic on $U-\{P\}$.
Then $z^n$ has a single-valued holomorphic logarithm on $U-\{P\}$,
which is absurd.
\end{proof}
\end{lemma}

\begin{corollary} $A$ and $B$ vanish identically.
\noproof
\end{corollary}

In fact, we do not exploit the vanishing of $B$, 
although it is a key step in the argument of \cite{G-SM} 
involving the linear system $\G_{00}$ of second
order theta functions that vanish to order $4$ at the origin 
and the heat equation.

Now $A$ can also be written as
\begin{eqnarray*}
A=&\frac{\partial}{\partial t}\bigg\vert_{t=0}\left(
\sum_{X\in S_g}a(X)N_g(\{x_{pq}\})\exp\pi i\sum_{pq,=1}^g x_{pq}(\tau_{pq}+
t\sigma_{pq})\right)\\
&=\frac{\partial}{\partial t}\bigg\vert_{t=0}N_g(F_g(\tau+t\sigma)).
\end{eqnarray*}
That is, $\sigma$ lies in the Zariski tangent space
$H$ at the point $\tau$ to the divisor in $\frak H_g$ defined
by the function $N_g(F_g)=N_g(\{{\frac{\partial}{\partial\tau_{ij}}}\})(F_g)$.
It is important to note that, from this description, $H$
depends upon $C$ but
is independent of the points $a,b$, the local co-ordinates
$z_a,z_b$ and the scalars $\lambda,\mu$.

We let $M_g^0$ denote the open subvariety of $M_g$
corresponding to curves with no automorphisms
and $A_g^0$ the open subvariety of $A_g$ corresponding to
principally polarized abelian varieties with no automorphisms
except $\pm 1$. Then $M_g^0$ lies in $A_g^0$ and both are
smooth varieties, and 
if $C$ lies in $M_g^0$ there are natural
identifications of tangent spaces given by
$$T_{[C]}M_g=H^0({\Omega^1_{C}}^{\otimes 2})^\vee,$$
$$T_{[C]}A_g =T_\tau\frak H_g=\Symm^2 H^0({\Omega^1_{C}})^\vee$$
(this latter identification is also a consequence of the heat equation).
The inclusion $T_{[C]}M_g\inj T_{[C]}A_g$
is dual to the natural multiplication (which is surjective,
by Max Noether's theorem)
$\Symm^2 H^0({\Omega^1_{C}})\to H^0({\Omega^1_{C}}^{\otimes 2}).$

We are aiming to prove that $H$,
when regarded as a Zariski tangent space,
is the whole of the tangent space 
$T_\tau\frak H_g=\Symm^2H^0(C,\Omega^1)^\vee.$
So assume otherwise; then $H$ is a hyperplane.
Projectivize: then 
$\sigma\in\P(H)$ and $\P(H)$ is a hyperplane
in $\P(\Symm^2H^0(C,\Omega^1)^\vee)$.

Now comes the point at which information about abelian integrals
is transformed into projective geometry and thence moduli.

The symmetric square $\Symm^2C$ is embedded in
$\P(\Symm^2 H^0(C,\Omega^1_C)^\vee)$
via the identification
$\Symm^2 H^0(C,\Omega^1_C)=H^0(\Symm^2C,{\Omega^1_C}^{\boxtimes 2})$,
where, by abuse of notation, 
${\Omega^1_C}^{\boxtimes 2}$ denotes the line bundle
on $\Symm^2C$
obtained by symmetrizing the exterior tensor square 
${\Omega^1_C}^{\boxtimes 2}$ on $C\times C$.
The entries $\sigma_{pq}$ of the matrix $\sigma$ 
are obtained by taking a basis of
$H^0(\Symm^2C,{\Omega^1_C}^{\boxtimes 2})$
and evaluating at the point $\{a,b\}$
of $\Symm^2C$. It follows,
since $H$ is independent of the points
$a$ and $b$, that the putative
hyperplane $\P(H)$ contains the embedded
$\Symm^2C$. However, $\Symm^2C$ is non-degenerate
in $\P(\Symm^2H^0(C,\Omega^1)^\vee)$
and therefore $H$ does not exist.
\end{proof}
\end{theorem}

Theorem \ref{non-transverse} is an immediate corollary
of this and the following lemma in commutative algebra.

\begin{lemma} \label{primary}
Suppose that $X$ is a closed subvariety of the variety $Y$
defined by the ideal $I=I_{X/Y}$. Suppose that $W$
is a smooth open subvariety of $Y$ such that $W\cap X$
is smooth and non-empty and that $J$ is an ideal
of $\sO_Y$ such that $J\vert_W=I^n\vert_W$.
Then $J$ is contained in $I^{[n]}$, the $n$th symbolic power of $I$.
\begin{proof} First, recall that if $X$ and $Y$ are
smooth over a field of characteristic zero,
then $I^n=I^{[n]}$ and consists of the functions
$f$ on $Y$ all of whose derivatives, with respect
to local co-ordinates on $Y$, of order up to
and including the $(n-1)$st, vanish along $X$.

We can assume that $Y$ is affine, say
$Y=\Sp A$, so that $A$ is an integral domain
and $I$ is prime. For any ideal $\frak a$ of
$A$, write $V(\frak a)=\Sp (A/\frak a)$. 

We can increase $J$, provided that $J\vert_W$ is unchanged,
so that in particular
we can replace $J$ by $J+I^{[n]}$. Then, without loss
of generality, we can suppose that $J$ contains $I^{[n]}$
and must prove that $J=I^{[n]}$. We have
$V(J)_{red}\subset V(I^{[n]})_{red}=X$ and
$V(J)_{red}\cap W=X\cap W$, so that
$V(J)_{red}=X$, and therefore ${\sqrt{J}}=I$.

Recall that for any ideal $\frak a$ with ${\sqrt{\frak a}}=I$, there
is a unique smallest $I$-primary ideal $\tfra$ containing $\frak a$,
given by the formula $\tfra =A\cap\frak a.A_I$, where
$A_I$ is the localization of $A$ at the prime ideal $I$. As before, we
can increase $J$, and so assume that $J={\widetilde{J}}$, that is,
that $J$ is $I$-primary. The symbolic power $I^{[n]}$
is $I^{[n]}={\widetilde{I^n}}$.

By assumption, the generic point $\xi$ of $X$ lies in $W$ and
$A_I=\sO_{Y,\xi}$, so that $J.A_I=I^n.A_I$. 
Intersecting both sides of this equation with $A$ gives
$J={\widetilde{J}}=I^{[n]}$.
\end{proof}
\end{lemma}

Now regard the Satake compactifications
$A_g^S$ and $M_{g+m}^S$ as
closed subvarieties of $A_{g+m}^S$.

\begin{theorem}\label{inf} (= Theorem \ref{non-transverse})
The intersection $A_g^S\cap M_{g+m}^S$
contains the $m$th order infinitesimal
neighbourhood of $M_g^S$ in $A_g^S$.
\begin{proof} The ideal defining $M_{g+m}^S$
inside $A_{g+m}^S$ is generated by those Siegel
modular forms $F_{g+m}$ that vanish along
$M_{g+m}^S$. From Theorem \ref{vanishing} and induction on $m$
it follows that $F_g$ and all its partial derivatives 
with respect to the co-ordinates $\tau_{pq}$ on
$\frak H_g$
of orders at most $m$ vanish along $M_g$, which is
just the statement of the corollary.
\end{proof}
\end{theorem}

\begin{remark} For $m=1$ this says
that at a general point $[C]$ of $M_g$, the Zariski
tangent space at $[C]$ to the $3g$-dimensional variety
$M_{g+1}^S$ contains the $g(g+1)/2$-dimensional tangent
space $\Symm^2H^0(C,\Omega^1_C)^\vee$ at $[C]$ to $A_g$,
where these tangent spaces both lie in $T_{[C]}A_{g+1}^S$.
\end{remark}

\end{section}
\bigskip

We are very grateful to Grushevsky and Salvati Manni for 
their correspondence that led to this paper and for their
interest in it, to Arbarello for some valuable conversations
and to an anonymous contributor to MathOverflow for the reference
to Bourbaki.
\bibliography{alggeom,ekedahl}
\bibliographystyle{pretex}
\end{document}